\documentclass{article}
\usepackage{a4,a4wide}
\usepackage{latexsym,amssymb,amsmath,amsthm}
\newtheorem*{theorem}{Theorem}
\newtheorem*{definition}{Definitions}
\newtheorem*{corollary}{Corollary}
\begin{document}
\author{Riccardo Poli and W. B. Langdon\\
Department of Computer Science\\
University of Essex, \\
Wivenhoe Park, Colchester, CO4 7SQ, UK\\
\{rpoli,wlangdon\}@essex.ac.uk}
\title{A simple reformulation of Riemann's Zeta function}
\date{}

\maketitle
\begin{abstract}
  We rewrite Riemann's Zeta function as 
  \[
  \zeta(s)=1+\sum_{p} p^{-s} \prod_{q\ge p} \frac{1}{1-q^{-s}}
\]
where $p$ and $q$ are primes.
\end{abstract}

\section{Introduction}

The connection between Riemann's zeta function $\displaystyle\zeta(s)=\sum_{n=1}^\infty n^{-s}$
\cite{riemann59:_ueber_anzah_primz_groes} and Euler's product
formula~\cite{euler48:_introd_analy_infin}  has long been known. In
this paper we introduce an intermediate representation for $\zeta(s)$
as a sum of products over the primes.

\section{Results}

\begin{definition}
  Let $\{p_k\}$ be the sequence of the primes,
\[
\mathcal Z_i(s)=\prod_{k=1}^i ( 1-p_k^{-s} )^{-1},
\]
\[
\mathcal S_i(s)=\sum_{k=1}^i \left ( p_k^{-s} \prod_{j=k}^i (1-p_j^{-s})^{-1} \right )
\]
and
\[
\mathcal I_i=  \{ s\in \mathbb C \,   | \, 
 p^{-s}=1 \mbox{ for some prime } p\le p_i \}.
\]
\end{definition}
Explicitly 
\[\mathcal I_i= \left \{ s\in \mathbb C \,  \Big | \, 
\Re(s)=0 \land \left ( \exists  k\in
\mathbb Z, j\in \{1,\cdots, i\}\, :\, \Im(s) = \frac{(1 + 2
  k)\pi}{\ln p_j}\right ) \right \}.
\]
\begin{theorem} For any $i \in \mathbb{N}$ and $s\in \mathbb
  C\setminus \mathcal I_i$, we have $Z_i(s) = 1 + S_i(s)$.
\end{theorem}
\noindent {\bf Proof:}
We proceed by induction. 

\noindent {\em Base case:} For  $i=1$ we have
\[
\mathcal S_1(s)=  p_1^{-s} ( 1-p_1^{-s} )^{-1} \Rightarrow
\]
\[
1+\mathcal S_1(s)= 1 + p_1^{-s} ( 1-p_1^{-s} )^{-1} =
\frac{1-p_1^{-s}+p_1^{-s}}{1-p_1^{-s}}  = \frac{1}{1-p_1^{-s}} =
\mathcal Z_1(s).
\]

\noindent {\em Induction hypothesis:} Assume 
\[
\mathcal{Z}_i(s) = 1 +\mathcal{S}_i(s).
\]
By adding and subtracting $p_{i+1}^{-s}$ to the r.h.s.\ we obtain
\begin{eqnarray*}
\mathcal{Z}_i (s)& = &   \left ( p_{i+1}^{-s} + \mathcal{S}_i(s) \right ) +
\left ( 1   - p_{i+1}^{-s} \right )
\Rightarrow \\[4mm]
(1   - p_{i+1}^{-s})^{-1} \mathcal{Z}_i(s)  &=&   (1   - p_{i+1}^{-s})^{-1}
\left ( p_{i+1}^{-s} + \mathcal{S}_i(s) \right )  +  (1   - p_{i+1}^{-s})^{-1} (1
- p_{i+1}^{-s} ) \Rightarrow 
\end{eqnarray*}
\begin{eqnarray*}
\mathcal{Z}_{i+1}(s)  &=&   (1   - p_{i+1}^{-s})^{-1}
\left ( p_{i+1}^{-s} + \sum_{k=1}^i \left ( p_k^{-s}  \prod_{j=k}^i (1-p_j^{-s})^{-1}\right ) \right )  +  1 \\
  &=&   (1   - p_{i+1}^{-s})^{-1}
 p_{i+1}^{-s} + \sum_{k=1}^i \left ( p_k^{-s}  \prod_{j=k}^{i+1} (1-p_j^{-s})^{-1}\right )   +  1 \\
  &=&    \sum_{k=1}^{i+1} \left ( p_k^{-s}  \prod_{j=k}^{i+1}
    (1-p_j^{-s})^{-1}  \right ) +  1 \Rightarrow 
\end{eqnarray*}

\[
\mathcal{Z}_{i+1}(s) = 1 +\mathcal{S}_{i+1}(s).
\]
\hfill$\Box$

\begin{corollary}
  For $\Re(s)>1$, Riemann's Zeta
  function can be written as
\[
\zeta(s)=1+\sum_{k=1}^\infty a_k(s)  p_k^{-s}
\]
where $\displaystyle a_k(s)=\prod_{j=k}^\infty (1-p_j^{-s})^{-1}$.
\end{corollary}

\noindent {\bf Proof:}
The result follows directly from  noting that Euler's product
formula  can be expressed as
\[
\zeta(s)=\lim_{i\rightarrow \infty}\mathcal Z_i
\]
and from  the previous theorem.
\hfill$\Box$

\section{Conclusions}

In our  reformulation for the Riemann's Zeta function, $\zeta(s)$ is
expressed as a sum over the primes, where each term is a product with
factors depending only on the summation index (a prime) and the primes
following it.  
The formulation is interesting because of its
similarity with the original definition of the Zeta function,
$\displaystyle\zeta(s)=\sum_{n=1}^\infty n^{-s}$, the factors 
$a_k(s)$ effectively being
corrections that account for all the (non-prime) naturals whose prime
factor decomposition starts with $p_k$.

\section*{Acknowledgements} The authors would like to thank John Ford
and Ray Turner for helpful discussions.


\end{document}